\documentclass[reqno,12pt]{amsart}
\usepackage{amsmath}
\usepackage{amscd}
\usepackage{amssymb}
\usepackage{enumerate}
\usepackage{amsfonts}
\usepackage{graphicx}
\usepackage[all]{xy}
\usepackage{mathrsfs}

\usepackage[mathscr]{euscript}
\usepackage{xypic}

\newtheorem{thm}{Theorem}
\newtheorem{prop}[thm]{Proposition}

\newtheorem{defn}[thm]{Definition}

\newcommand{\bC}{\mathbb{C}}

\newcommand{\bP}{\mathbb{P}}

\newcommand{\bZ}{\mathbb{Z}}

\newcommand{\sN}{\mathscr{N}}
\newcommand{\sO}{\mathscr{O}}
\newcommand{\sP}{\mathscr{P}}

\newcommand{\sV}{\mathscr{V}}

\newcommand{\sX}{\mathscr{X}}

\newcommand{\fb}{\mathfrak{b}}
\newcommand{\fsl}{\mathfrak{sl}}
\newcommand{\fh}{\mathfrak{h}}
\newcommand{\fg}{\mathfrak{g}}

\newcommand{\KS}{\chi}

\newcommand{\bung}{{{\mathscr M}_G}}
\newcommand{\cbg}{\sX}

\DeclareMathOperator{\Grass}{Grass}

\DeclareMathOperator{\Prym}{{Prym}}
\DeclareMathOperator{\rank}{rank}
\DeclareMathOperator{\Res} {Res}

\DeclareMathOperator{\Sym} {Sym}
\DeclareMathOperator{\Tr}  {Tr}

\title{Donagi-Markman cubic for Hitchin systems}
\bibliography{ref}

\author[D.~Balduzzi]{David Balduzzi}

\thanks{Partially supported by NSF grant DMS-0401164. \\
{\em Address:} Department of Mathematics, University of Chicago, Chicago, IL 60637 \\
{\em E-mail:} {\tt balduzzi@math.uchicago.edu}}

\begin{document}

\maketitle

\begin{abstract}
The Donagi-Markman cubic is the differential of the period map for algebraic completely integrable systems. Here we prove a formula for the cubic in the case of Hitchin's system for arbitrary semisimple $\fg$. This was originally stated (without proof) by Pantev for $\fsl_n$.
\end{abstract}

\section{Introduction}
\label{sec:intro}

Hitchin's system is an integrable system constructed out of the moduli space of principal bundles on a Riemann surface. The cotangent bundle to the moduli space is canonically a symplectic manifold, points of which are principal bundles together with a twisted section known as the Higgs field. A Lagrangian fibration is found by mapping such a pair to the ``characteristic polynomial'' of the Higgs field. Fibers of this map are open subsets of Prym varieties of a cover of the Riemann surface. This is described in more detail in \S\ref{sec:spec}.

Given a family of Abelian varieties, one can try calculate their periods. For the Hitchin system this seems too hard. However the differential of the period map is considerably easier to find.
For any lagrangian fibration by abelian varieties it turns out that this differential forms a cubic tensor on the base, the Donagi-Markman cubic \S\ref{sec:cub}.

In unpublished notes T. Pantev gave (without proof) a formula for the cubic for ${\mathfrak sl}_n$. Our goal in this note is to prove the formula in the general case. The proof is essentially a Kodaira-Spencer map computation. To handle arbitrary semisimple groups we work with cameral covers rather than the more familiar spectral covers.

Thus our goal is to prove the following

\begin{thm}
\label{thm:pantev}
The cubic map (\ref{eq:cubhit}) for the Hitchin system
\[
c_\phi:H^0(\tilde{M}_\phi, \fh\otimes K)^W\rightarrow \Sym^2 \left[ H^0(\tilde{M}_\phi, \fh\otimes K)^W\right]^*
\]
is given by Pantev's formula (\ref{eq:pantev3})
\begin{equation*}
\beta\mapsto\left(\gamma\cdot\delta \mapsto\Res^2_{\tilde{M}_\phi}\left[\pi^*\left( \frac{d_\phi D(\beta)}{D(\phi)}\right)\Tr\left(\gamma\cup\delta\right)\right]\right),
\end{equation*}
for $\beta, \gamma$ and $\delta\in H^0(\tilde{M}_\phi,\fh\otimes K)^W$.
\end{thm}

\S\ref{sec:pantev} explains the formula. The proof is contained in \S\ref{sec:bpc} and \S\ref{sec:proof}. Finally in Theorem \ref{thm:symm} we produce a more symmetric formula for the cubic:
\begin{equation*}
c_\phi:(\beta,\gamma,\delta)\mapsto\sum_{\nu\in\Delta}
\Res_{\tilde{M}_\phi}\left[\frac{<\beta,\nu><\gamma,\nu><\delta,\nu>}{<\phi,\nu>}\right]
\end{equation*}

{\bf Acknowledgements.}
I'd like to thank my advisor, V. Ginzburg, for suggesting this problem and for many helpful conversations.

\section{Hitchin's system}
\label{sec:spec}

Given an arbitrary complex connected semisimple Lie group $G$ with Lie algebra $\fg$, let $\bC[\fg]$ be the algebra of polynomials on $\fg$. Choosing a Cartan subalgebra $\fh\subset \fg$ we can form the restriction map $\bC[\fg]\rightarrow \bC[\fh]$. There is a natural $G$-action on $\fg$ given by the adjoint representation. The Cartan $\fh$ comes with a natural action of the Weyl group $W$ of $G$. We can then consider $G$ and $W$-invariant polynomials on the respective spaces. Chevalley's theorem states that restriction induces an isomorphism $\iota:\bC[\fg]^G\xrightarrow{\approx}\bC[\fh]^W$.
In addition it is known that $\bC[\fg]^G$ is a free polynomial algebra, generated by homogeneous polynomials $p_i$ of degree $d_i\geq 2$ for $i=1,\ldots,r=\rank\fg$.

Suppose $M$ is a compact Riemann surface of genus $>1$ and $G$ a semisimple Lie group. Consider the moduli space of stable principal $G$-bundles of degree $d$ on $M$. This is generically a smooth manifold with finite quotient singularities at principal bundles where the center $Z(G)$ does not contain all automorphisms of the bundle, see \cite{ram_curveI} and \cite{ram_curveII}.  We denote the smooth locus by $\bung$, with tangent space at a principal bundle $P$ given by $H^1(M,\fg_P)$, where $\fg_P$ is the vector bundle associated to $P$ via the adjoint representation.
By Serre duality, the cotangent space is $H^0(M, \fg_P\otimes K_M)$ where $K_M$ is the canonical line bundle on $M$.

Picking (non-canonically) generators $p_i$ for $\bC[\fg]^G$, define a twisted version of the Chevalley map
\begin{align}
\label{eq:h}
h_P:H^0(M,\fg_P\otimes K_M)&\rightarrow \bigoplus_{i=1}^r H^0(M,K_M^{d_i}) =: B\\
\Phi&\mapsto \sum_i  p_i(\Phi)=:\sum_i \phi_i,\,\,\mbox{ where }\phi_i\in H^0(M,K_M^{d_i})
\end{align}
and so obtain the Hitchin map $h:T^*\bung\rightarrow B$. A {\it Higgs bundle} is defined to be a pair $(P,\Phi)\in T^*\bung=:\cbg$. Both $\Phi$ and $\phi$ are referred to as the {\it Higgs field}.  Since $\cbg$ is a cotangent space to a manifold it is has a natural symplectic structure and Hitchin shows \cite{hit_is} that $h$ is a Lagrangian fibration. The tangent space to $\cbg$ at a point $(P,\Phi)$ is described by exact sequence
\[
0\rightarrow H^0(M, \fg_P\otimes K_M)\rightarrow T_{(P,\Phi)}\cbg\rightarrow H^1(M,\fg_P)\rightarrow 0.
\]

Next we describe the fibers of $h$. In the case of the classical groups, for each element of the base $b\in B$ Hitchin constructs a cover of $M$. Line bundles on this cover correspond (roughly) to $G$-bundles on $M$, and some subset of the Picard group of this cover gives the fiber of $h$ over $b$. For arbitrary semisimple $G$, it is useful to follow Donagi and Scognamillo (see \cite{d_spec}, \cite{d_decomp} and \cite{sc_abel}) and use cameral covers.

Fix Cartan $H$ and Borel $B$ in $G$. Consider again Chevalley's restriction map, this time from a scheme-theoretic perspective where it is known as the adjoint quotient map. It provides a morphism of affine varieties
\[
\fg\twoheadrightarrow \fh/W
\]
and from this we can form the fiber product
\[
\xymatrix{
\fg\times_{\fh/W}\fh=\tilde{\fg}\ar[r]\ar[d] & \fh \ar[d] \\
\fg \ar[r] & \fh/W.
}\]

We wish to globalize this construction to vector bundles on $M$. Let $|\cdot|$ denote the total space of a vector bundle on $M$. Twisting by the canonical bundle of $M$ gives a relative version of the diagram
\[
\xymatrix{
\left|\fg_P\otimes K\right|\times_{\left|\fh\otimes K\right|/W}\left|\fh\otimes K\right|
\ar[r]\ar[d] & \left|\fh\otimes K_M\right| \ar[d] \\
\left|\fg_P\otimes K_M\right| \ar[r] & \left|\fh\otimes K_M\right|/W
}\]
and the cameral cover of $M$ is defined as the pullback
\[
\tilde{M}:= \Phi^*\left(\left|\fg_P\otimes K\right|\times_{\left|\fh\otimes K\right|/W}
\left|\fh\otimes K\right|\right).
\]

Alternatively by Chevalley's theorem $\tilde{M}=\phi^*\left(\left|\fh\otimes K_M\right|\right)$. From this it follows that $\tilde{M}$ comes with a natural $W$-action and a $W$-equivariant $M$-morphism $\tilde{M}\rightarrow \left|\fh\otimes K\right|$.
Projecting on either side of the commutative diagram gives a map $\tilde{M}\xrightarrow{\pi}M$, since $\Phi$ or $\phi$ are sections of $\fg_P\otimes K$ or $\left|\fh\otimes K\right|/W$ respectively. Over the points $m$ where $\Phi$ is regular semisimple the fiber $\pi^{-1}(m)$ is the set of chambers in $\fh^*_m$, hence the name cameral cover. In particular for ${\mathfrak sl}_n$, the fiber is given by the set of orderings of the eigenvalues of $\Phi(m)$. By abuse of notation we use $K_M$ to denote the line-bundle $\pi^* K_M$ on $\tilde{M}$.

The key fact about the cameral curve is that it {\it abelianizes} the Higgs bundle $(P,\Phi)$.
The abelianization procedure can be cast in the abstract setting of principal Higgs bundles, see \cite{dg_gerbe}. However since we are considering concrete $K$-valued Higgs bundles it suffices to define cameral covers as a globalized form of the $W$-cover $\fh\rightarrow \fh/W$.
 
The pullback $\pi^*P$ has natural reduction to principal bundle $P_B$ with structure group $B$ so that $\pi^*\Phi\in H^0(\tilde{M},\fb_{P_B}\otimes K_M)$, see \S3 in \cite{sc_abel}. The projection map $B\rightarrow H$ associates to $P_B$ principal $H$-bundle $P_H$, and this can be twisted (see \cite{sc_abel}) to $W$-equivariant $H$-bundle $\tilde{P}_H$. So we have a map
\begin{equation}
\label{eq:lociso}
(P,\Phi)\mapsto(\tilde{M},\tilde{P}_H).
\end{equation}

\begin{thm}
\cite{hm_pr}
\renewcommand{\labelenumi}{\alph{enumi})}
\begin{enumerate}
\item
Let $\tilde{M}_0$ be cameral cover. Let $\sN$ be the variety of pairs $(\tilde{M}, \tilde{P}_H)$ where $\tilde{M}$ is a $W$-invariant deformation of $\tilde{M}_0$ in $\left|\fh\otimes K\right|$, and $\tilde{P}_H$ is a $W$-invariant $H$-bundle on $\tilde{M}$. Then under the map (\ref{eq:lociso}) $\sN$ is locally isomorphic to $\cbg$.
\item
The projection $(\tilde{M},\tilde{P}_H)\mapsto \tilde{M}$ defines a Lagrangian fibration of open subset of $\cbg$.
\end{enumerate}
\end{thm}

Let $\Prym(\tilde{M})$ be the variety of $W$-invariant $H$-bundles on $\tilde{M}$. The fiber $h^{-1}(b)$ of the Hitchin map is an open subset of $\Prym(\tilde{M}_b)$, but note this subset does {\it not} lie in the connected component containing the identity. The Lagrangian foliation gives rise to exact sequence
\[
0\rightarrow H^1(\tilde{M},\sO\otimes\fh)^W\rightarrow T\cbg\rightarrow H^0(\tilde{M},\fh\otimes K_M)^W\rightarrow 0,
\]
where $H^1(\tilde{M},\fh \otimes\sO)^W$ is the tangent space to $\Prym$ and $H^0(\tilde{M},\fh\otimes K_M)^W$ is the tangent space to $B$. The embedding of $\tilde{M}$ in $\left|\fh\otimes K_M\right|$ produces exact sequence
\begin{equation}
\label{eq:embcam}
0\rightarrow T\tilde{M}\rightarrow T\left|\fh\otimes K_M\right|\rightarrow N_{\tilde{M}}\rightarrow0
\end{equation}
Infinitesimal $W$-deformations of a cameral curve $\tilde{M}$ embedded in $\left|\fh\otimes K\right|$ correspond to sections of $H^0(\tilde{M},N_{\tilde{M}})^W$.
Thus if we consider $B$ to be the moduli space of cameral curves, it follows that $T_b B=H^0(\tilde{M}_b,N)^W$.

There is a natural $\fh$-valued two-form on $\left|\fh\otimes K_M\right|$ induced by the symplectic structure on the total space $|K_M|$. This gives rise to map $N_{\tilde{M}}\rightarrow \fh\otimes K_M$ which induces isomorphism $H^0(\tilde{M},N_{\tilde{M}})^W\rightarrow H^0(\tilde{M},\fh\otimes K_M)^W$,  see \cite{hm_pr} . In the sequel we will use this identification $T_bB=H^0(\tilde{M}_b,\fh\otimes K_M)^W$.

\section{The cubic}
\label{sec:cub}

Now we take a look at Hitchin's system from a more general point of view. Let $B^{sm}\subset B$ be the open subset of $B$ corresponding to smooth cameral curves where the Higgs field has at most simple zeroes. For the rest of this section work locally on $B$ in a neighborhood of a base point $0\in B^{sm}$. Abusing notation we tend to write $B$ for this neighborhood. Denote the fiber $\pi^{-1}(b)$ by $X_b$. Since we are restricting to the smooth locus, and the Hitchin map $h$ is a proper submersion there is a diffeomorphism
\[
T:\sX\approx X\times B
\]
where $X:=X_0$. In particular, all the fibers are diffeomorphic. Thus we can speak about varying the holomorphic structure on the base fiber $X$ instead of varying the fiber itself. There is an exact sequence on $\cbg$ relating the base and relative tangent bundles
\begin{equation}
\label{eq:reltan}
0\rightarrow T_{\sX/B}\rightarrow T_\sX\rightarrow \pi^*T_B\rightarrow 0.
\end{equation}

\begin{defn}
\label{df:ks}
Restricting the resulting sheaf cohomology long exact sequence to the fiber $X$, there is coboundary map
\begin{equation*}
\KS:T_{B,0}=H^0(X,\pi^* T_{B|X})\rightarrow H^1(X,T_X).
\end{equation*}
$\KS$ is the Kodaira-Spencer map, measuring the infinitesimal variation of the holomorphic structure on $X$.
\end{defn}

Another way to study the variation of the holomorphic structure is to use the Hodge filtration. Consider the inclusion $H^0(X_b,\Omega^1)\subset H^1(X_b,\bC)$. Here $H^1(X_b,\bC)$ depends only on the topology of $X_b$, and since $X_b$ is diffeomorphic to $X$, we can identify it with $H^1(X,\bC)$ via the Gauss-Manin connection.

\begin{defn}
\label{df:period}
(Griffiths) The period map
\[
\sP:B\rightarrow\Grass(d,H^k(X,\bC)),
\]
where $d=\dim H^0(X_b,\Omega^1)$, is the map which to $b\in B$ associates the subspace
\[
H^0(X_b,\Omega^1)\subset H^1(X_b,\bC)\approx H^1(X,\bC).
\]
\end{defn}

Choosing bases of holomorphic differentials and 1-cycles, this can be shown to reduce to the usual period map on, for example, curves or abelian varieties. The period map and the Kodaira-Spencer map are related by the following

\begin{thm}
\label{thm:per}
(Griffiths, \cite{gr_p2} Theorem 1.27)
Given $u\in T_{B,0}$ the map
\[
d\sP(u):H^0(X,\Omega^1_X)\rightarrow H^1(X,\sO_X)
\]
is equal to the cup-product with the class $\KS(u)\in H^1(X,T_X)$, composed with the map induced on cohomology by the interior product $T_X\otimes\Omega^1_X\rightarrow \sO_X$.
\end{thm}

This map measures the variation in complex structure of the abelian variety\\ $A=H^1(X,\sO)/H^1(X,\bZ)$ associated to $X$ -- the variation in the periods. Finally Griffiths shows that if there is a family of polarizations $\left(\pi_* \Omega^1_{\sX/B}\right)^*\xrightarrow{\approx} R^1\pi_*\sO_\sX$ then Riemann's symmetric bilinear relation holds
\begin{equation}
\label{eq:r1}
d\sP:T_{B,0}\rightarrow \Sym^2 H^1(X,\sO_X)\subset \otimes^2 H^1(X,\sO_X)
\end{equation}
after using the identification given by the polarization.

In \cite{dm_cub} Donagi and Markman investigate conditions on a family of abelian varieties $\pi:\cbg\rightarrow B$ so that locally on $B$ there is a symplectic structure on $\cbg$ such that $\pi$ is a Lagrangian fibration. Let $\sV$ be the bundle $\pi_* T_{\cbg/B}$ on $B$ with fibers the tangent space to the fibers of $\pi$. Since the fibers are Lagrangian subvarieties,  the symplectic form $\omega$ on $\cbg$ induces an isomorphism
\[
\iota:\sV^*\rightarrow TB
\]
The cubic is the composition $c:=d\sP\circ\iota:\sV^*\rightarrow \Sym^2 \sV$, and in the case of a Lagrangian fibration, it turns out by a Koszul-complex computation that
\[
c\in H^0(B, \Sym^3 \sV) \subset H^0\left(B,\sV\otimes \Sym^2 \sV\right).
\]

Our aim is to calculate the cubic in the case of Hitchin's system.  Let $P_b$ be the Prym variety corresponding to a point $b\in B$. This inherits polarization (not in general principal) by Serre duality on $\tilde{M}$ since $\left[H^0(\tilde{M}_b,\fh\otimes K)^W\right]^*\approx TP_b$. The cubic goes
\begin{equation}
\label{eq:pantev1}
c_b:T_b B\rightarrow H^0(P,\Sym^2 TP),
\end{equation}
so identifying $T_b B\approx H^0(\tilde{M}_b,\fh\otimes K)^W$ as in the previous section the cubic is a map
\[
c_b:H^0(\tilde{M}_b,\fh\otimes K)^W\rightarrow \Sym^2 H^1(\tilde{M}_b,\fh\otimes\sO)^W.
\]

It is constructed as follows. The Kodaira-Spencer map
\begin{equation}
\label{eq:KSb}
\KS:T_b B=H^0(\tilde{M}_b,N)^W=H^0(\tilde{M}_b,\fh\otimes K)^W\rightarrow H^1(\tilde{M}_b,T)
\end{equation}
measures deformations of the cameral cover. We are interested in deformations of the associated Prym. 
To this end consider the following map constructed out of the Killing form on $\fg$, which we denote by $\Tr$, and the cup product:
\begin{equation}
\label{eq:cuptr}
\Sym^2 H^0(\tilde{M},\fh\otimes K)^W\rightarrow \otimes^2 H^0(\tilde{M},\fh\otimes K)^W\xrightarrow{\cup} H^0(\tilde{M},\fh\otimes \fh\otimes K^2_M)\xrightarrow{\Tr} H^0(\tilde{M},K_M^2).
\end{equation}
Denote the (Serre) dual map by $G$:
\begin{equation}
\label{eq:kp}
G:H^1(\tilde{M},T)\rightarrow\otimes^2 H^1(\tilde{M},\fh\otimes \sO)^W\rightarrow\Sym^2 H^1(\tilde{M},\fh\otimes \sO)^W.
\end{equation}
Theorem \ref{thm:per} can be rewritten in this situation as
\begin{prop}
\label{prop:G}
Given $\beta\in T_{B,0}$, the differential of the period map
\[
d\sP(\beta):H^0(\tilde{M},\fh\otimes K)^W\rightarrow H^1(\tilde{M},\fh\otimes\sO)^W
\]
measuring variation of holomorphic structure on the Prym is given by cup-product with the class $\chi(\beta)\in H^1(\tilde{M},T)$, composed with $G$:
\begin{equation*}
d\sP:T_{B,0} \rightarrow \Sym^2 H^1(\tilde{M},\fh\otimes\sO)^W : \beta\mapsto G\circ\chi(\beta).
\end{equation*}
\end{prop}


Recalling the identification of $T_b B$ with $H^0(\tilde{M}_b, \fh\otimes K)^W$ we obtain the cubic. This can be rewritten in a form more similar to (\ref{eq:pantev3}) using (\ref{eq:cuptr})

\begin{align}
\label{eq:cubhit}
c_b:H^0(\tilde{M}_b,\fh\otimes K)^W & \rightarrow\Sym^2\left[H^0(\tilde{M}_b,\fh\otimes K)^W\right]^* \\
\beta & \mapsto\left(\gamma\cdot\delta \mapsto \int_{\tilde{M}_b}\left[\KS(\beta)\cup \Tr(\gamma\cup\delta)\right]\right),
\end{align}
where here $\int:H^1(\tilde{M}_b,K_{\tilde{M}})\rightarrow \bC$.

\section{Branch points and the discriminant}
\label{sec:pantev}

Given a semi-simple Lie algebra $\fg$ and Cartan subalgebra $\fh$ let $\Delta$ be the set of roots of $\fg$. Define the discriminant as
\begin{equation}
\label{eq:disc}
D:\fh\rightarrow\bC:h\mapsto\prod_{\gamma\in\Delta}\gamma(h).
\end{equation}
Since we multiply over all roots in $\fh$ the discriminant is independent of any choices. It cuts out a divisor in $\fh$, also denoted by $D$. For ${\mathfrak sl}_n$ if we think of $\fh$ as the collection of ordered $n$-tuples of eigenvalues, then $D$ is the divisor where two or more of the eigenvalues in a tuple coincide. Given a Higgs bundle $(P,\Phi)$, globalize to
\[
M\xrightarrow{h_P(\Phi)}\left|\fh\otimes K\right|/W\xrightarrow{D}|K^{|\Delta|}|,
\]
where $|\Delta|$ is the number of roots. Pulling back the zero-section of $K^{|\Delta|}$ cuts out a divisor on $M$. This divisor marks where two or more eigenvalues of the Higgs field coincide: branch-points of the cameral cover $\tilde{M}$.
Varying the Higgs field causes corresponding variation in the cameral cover and this variation is encoded in the motion of the branch points. Extending the base field by a nilpotent $\epsilon$ satisfying $\epsilon^2=0$ and picking $\beta$ in $T_\phi B\approx B$ we can deform the Higgs field:
\[
\phi+\epsilon\cdot\beta:M\rightarrow \fh/W\otimes K \xrightarrow{D}K_M^{|\Delta|}.
\]
This deformation gives a 1-parameter family of divisors on $M$, the family of branch-points coming from the family of cameral covers. For ${\mathfrak sl}_n$ the divisor looks like
\begin{equation}
\label{eq:br_div}
\prod_{i\neq j}\left(\phi_i-\phi_j+\epsilon(\beta_i-\beta_j)\right) =D(\phi)+\epsilon\cdot d_\phi D(\beta)+\mbox{ higher order terms }\ldots
\end{equation}
where the higher order terms vanish since $\epsilon^2=0$.

Note the terms on the left-hand side of this equation are only well-defined locally on the Hitchin base since there is a non-trivial monodromy group acting on the branch points. The product is well-defined globally since the divisor cut out by the discriminant in independent of the choice of labeling of the roots.

Thus we obtain a linear system of sections of $H^0(M,K^{|\Delta|})$ generated by $D(\phi)$ and $d_\phi D(\beta)$. This 1-dimensional linear system gives rise to a meromorphic function
\begin{equation}
\label{eq:discratmer}
\frac{d_\phi D(\beta)}{D(\phi)}:M\rightarrow\bP^1,
\end{equation}
the discriminant ratio. Using this construct map
\begin{align}
\label{eq:pantev2}
H^0(\tilde{M}_\phi,\fh\otimes K_M)^W&\rightarrow \Sym^2 \left[H^0(\tilde{M}_\phi, \fh\otimes K_M)^W\right] ^*
\\
\label{eq:pantev3}
\beta&\mapsto\left(\gamma\cdot\delta\mapsto\Res^2_{\tilde{M}_\phi}\left[\pi^*\left(\frac{d_\phi D(\beta)}{D(\phi)}\right)\Tr\left(\gamma\cup\delta\right)\right]\right).
\end{align}
where $\beta, \gamma$ and $\delta \in H^0(\tilde{M}_\phi, \fh\otimes K)^W$.
Here $\Res^2$ is the quadratic residue on a Riemann surface, locally described as
\begin{equation*}
\Res_{\{0\}}^2(\omega)=\alpha,\mbox{ for }
\omega=\alpha\frac{(dz)^2}{z^2}+\beta\frac{\left(dz\right)^2}{z}+\ldots
\end{equation*}

The discriminant ratio is pulled back from $M$ to $\tilde{M}_\phi$. It has simple poles on the branch points of $M$, and when pulled-back by $\pi$ these poles become quadratic since near branch points $\pi$ looks like $z\mapsto z^2$. This follows from the assumption that the Higgs field only hits the smooth points of the
discriminant divisor.

In the next two sections we look at deformations of a branched double cover and then prove Theorem \ref{thm:pantev}.

\section{A Kodaira-Spencer computation}
\label{sec:bpc}

The Kodaira-Spencer map measures the infinitesimal variation in holomorphic structure in a family of complex curves. We are studying a family of cameral curves over the base curve $M$, and the holomorphic structure on the cameral curve is determined by the pattern of branch points. So we begin by looking at the local picture determined by nudging a branch point.

The equation
\begin{equation}
\label{eq:eg2}
\lambda^2=q(z),\,\,\mbox{ where } q(0)=0,
\end{equation}
determines a two sheeted covering of the $z$-plane corresponding to the solutions $\pm\sqrt{q(z)}$. In addition assume $\dot{q}(0)\neq0$ so that the cover is smooth. Suppose we deform this equation as follows
\begin{equation}
\label{eq:eg2def}
\lambda^2=q(z)+\epsilon \beta(z),
\end{equation}
giving solutions
\begin{equation}
\label{eq:sol2}
\lambda=\pm\sqrt{q(z)+\epsilon\beta(z)}=\pm\left(\sqrt{q(z)}+\epsilon\frac{\beta}{2\sqrt{q(z)}}\right).
\end{equation}
How does the cover vary as a complex manifold? Here we introduce the Kodaira-Spencer map. In general suppose $z_i$'s are a collection of local coordinates on a Riemann surface with transition functions $f_{ij}$ depending on a parameter $\epsilon$.
\begin{gather*}
z_i=f_{ij}(z_j,\epsilon)\,\,\mbox{ and }\,\,z_j=f_{jk}(z_k,\epsilon),\\
\mbox{so that }\,\, f_{ik}(z_k,\epsilon)=f_{ij}(f_{jk}(z_k,\epsilon),\epsilon).
\end{gather*}
Then
\begin{align*}
\frac{\partial f_{ik}}{\partial \epsilon}(z_k,\epsilon)&=
\frac{\partial f_{ij}}{\partial \epsilon}(z_j,\epsilon)+\frac{\partial f_{ij}}{\partial z_j}(z_j,\epsilon)\cdot\frac{\partial f_{jk}}{\partial \epsilon}(z_k,\epsilon)\\
&=\frac{\partial f_{ij}}{\partial \epsilon}(z_j,\epsilon)+\frac{\partial z_i}{\partial z_j}\cdot\frac{\partial f_{jk}}{\partial \epsilon}(z_k,\epsilon).
\end{align*}
Since $\frac{\partial}{\partial z_j}=\frac{\partial z_i}{\partial z_j}\frac{\partial}{\partial z_i}$ this is equivalent to
\begin{equation}
\label{eq:coc_hid}
\frac{\partial f_{ik}}{\partial \epsilon}(z_k,\epsilon)\frac{\partial}{\partial z_i}=\frac{\partial f_{ij}}{\partial \epsilon}(z_j,\epsilon)\frac{\partial}{\partial z_i}+\frac{\partial f_{jk}}{\partial \epsilon}(z_k,\epsilon)\frac{\partial}{\partial z_j}.
\end{equation}
Introduce 1-cycles
\begin{equation}
\label{eq:1cyc}
\theta_{jk}:=\frac{\partial f_{jk}}{\partial \epsilon}(z_k,\epsilon)\frac{\partial}{\partial z_j}.
\end{equation}
Now notice that (\ref{eq:coc_hid}) is equivalent to the cocycle condition
\begin{equation}
\theta_{ik}=\theta_{ij}+\theta_{jk}.
\end{equation}
This give the Kodaira-Spencer map from the moduli space of deformations of a complex variety $S$, to $H^1(S,T_S)$. In the case (\ref{eq:eg2def}) above, we have
\begin{gather*}
\lambda=\pm\left(\sqrt{q}+\frac{\epsilon\beta}{2\sqrt{q}}\right)\,\,\mbox{ on the two branches,}\\
\mbox{so }\,\, \theta^1_{\lambda z}=\frac{\beta(z)}{2\sqrt{q(z)}}\frac{\partial}{\partial\lambda} = \frac{\beta}{2\lambda}\frac{\partial}{\partial\lambda}\,\,\mbox{ and similarly }\,\,\theta^2_{\lambda z}=\frac{\beta}{2\lambda}\frac{\partial}{\partial \lambda}.
\end{gather*}
Thus it follows that the Kodaira-Spencer cocycle is
\begin{equation}
\label{eq:theta_lz}
\theta_{\lambda z}=\frac{\beta}{\lambda}\frac{\partial}{\partial\lambda}.
\end{equation}
We want to link this 1-cocycle with the discriminant: $D(\lambda^2-q(z)-\epsilon\beta(z))=4q(z)+4\epsilon\beta(z)$, so
\begin{equation}
\label{eq:drat2}
\frac{d_\epsilon D(\beta)}{D(q)}=\frac{\beta(z)}{q(z)}=\frac{\pm\beta/\sqrt{q}}{\pm\sqrt{q}}=\frac{\beta}{\lambda^2}.
\end{equation}

Given a holomorphic quadratic differential $f(\lambda)(d\lambda)^2$ defined on the cover in a neighbourhood of the branch point, we can multiply it by the discriminant ratio (\ref{eq:drat2}) and take the quadratic residue. Alternatively we can evaluate it against the Kodaira-Spencer cocycle (\ref{eq:theta_lz}). This gives a meromorphic differential with a $1^{st}$-order pole at the branch point. Taking its residue gives the same answer:
\begin{equation}
\label{eq:compare}
\Res^2_{\lambda=0}\left(f(\lambda)\frac{\beta}{\lambda^2}\right)(d\lambda)^2 = \Res_{\lambda=0}\left(f(\lambda)\frac{\beta}{\lambda}\right)d\lambda.
\end{equation}

\section{Proof of the formula}
\label{sec:proof}

\begin{proof}

Given root $\nu\in\Delta\subset\fh^*$ and $\alpha\in\fh$, write $<\alpha,\nu>$ for the evaluation pairing. In terms of roots the meromorphic discriminant ratio $M\rightarrow\bP^1$ is then
\begin{align}
\label{eq:discrat}
\frac{d_\phi D(\beta)}{D(\phi)} = \frac{d_\phi \left(\prod_{\nu\in\Delta}<\beta,\nu>\right)}{\prod_{\nu\in\Delta}<\phi,\nu>}
& = \sum_{\nu\in\Delta}\frac{<\beta,\nu>}{<\phi,\nu>}, \\
\mbox{which for }{\mathfrak sl}_n\mbox{ is } & = \sum_{i\neq j}\frac{\beta_i-\beta_j}{\phi_i-\phi_j}.
\end{align}
The individual terms on the right-hand side only make sense after pulling back to the cameral cover, however the expression as a whole is well-defined on $M$. Working in a neighbourhood of a branch point where $<\alpha,\nu_1>=0$ -- equivalently near a pole where $\phi_1$ and $\phi_2$ collide --  write this as
\begin{equation}
\label{eq:discratloc}
\frac{\beta_1-\beta_2}{\phi_1-\phi_2}+\mbox{ regular terms }\ldots
\end{equation}
where ``regular terms'' are holomorphic near the branch point and so will not contribute.
All branch points are quadratic since we are assuming our Higgs fields to have at most simple zeroes. So we reduce the situation near a branch point to that of the calculation in \S\ref{sec:bpc}.
Near the branch point, choose a local $z$ coordinate so that $z=0$ is the branch point, and the equation defining the cameral cover looks like
\[
(\lambda^2 - q(z))(\mbox { regular terms })
\]
as in \S\ref{sec:bpc}. Then (\ref{eq:discratloc}) becomes
\begin{equation}
\label{eq:drpc}
\frac{\beta}{q}+\mbox{ regular terms }\ldots
\end{equation}
Equation (\ref{eq:theta_lz}) calculates the Kodaira-Spencer map giving
\begin{align*}
\theta_{\lambda z}&=\frac{\beta(z)}{\sqrt{q(z)}}\frac{\partial}{\partial\lambda},\,\,\mbox{ and so}\\
c_\phi:T_\phi B\rightarrow & H^1(\tilde{M}_\phi,T):
\beta\mapsto\frac{\beta_1-\beta_2}{\phi_1-\phi_2}\frac{\partial}{\partial \lambda}.
\end{align*}

Using (\ref{eq:compare}) we see that the quadratic residue of a quadratic differential times the discriminant ratio coincides with the normal residue of the quadratic residue cupped with the image of the Kodaira-Spencer. Standard argments show that this in turn coincides with the map $\Res:H^1(\tilde{M},K_{\tilde{M}})\rightarrow\bC$ given by integration.
\end{proof}

In the special case of $\fg=\fsl_2$ Pantev derives a simpler version of the formula which depends only on the base curve $M$:
\begin{align*}
H^0(M, K_M^2)&\rightarrow \Sym^2 H^0(M,K_M^2)\\
\beta&\mapsto \left(\gamma\mapsto\Res^2_M\left[\frac{\beta\cdot\gamma^2}{\phi^2}\right]\right).
\end{align*}
To recover this formula, we return to the calculations in \S\ref{sec:bpc}. We had deformation
$\lambda^2=\phi(z)+\epsilon\beta(z)$ giving rise to solutions
\[
\lambda=\pm\left(\sqrt{\phi}+\frac{\epsilon\beta}{2\sqrt{\phi}}\right).
\]
For ${\mathfrak sl}_2$ the spectral and cameral covers coincide and the cameral cover embeds in $|\fh\otimes K_M|=|K_M|$. Performing Lagrange interpolation, see \cite{hurt_is}, gives map
\[
K^2_M\rightarrow K_{\tilde{M}}:\beta\mapsto \frac{\beta}{\sqrt{\phi}}
\]
So when we calculate $\Res^2_{\tilde{M}}\left(\frac{d_\phi D(\beta)}{D(\phi)}\gamma^2\right)$ in terms of the base we are calculating
\[
\Res^2_M\left(\frac{\beta}{\phi}\left(\frac{\gamma}{\sqrt{\phi}}\right)^2\right).
\]
\newline

At first glance the cubic is not symmetric in $\beta, \gamma$ and $\delta$ -- it is not obviously a cubic. However we have the following reformulation of the formula

\begin{thm}
	\label{thm:symm}
	The cubic for the Hitchin system is given by formula
	\begin{equation*}
		c_\phi:(\beta,\gamma,\delta)\mapsto\sum_{\nu\in\Delta}
		\Res_{\tilde{M}_\phi}\left[ \frac{<\beta,\nu><\gamma,\nu> <\delta,\nu>}{<\phi,\nu>} \right].
	\end{equation*}	
\end{thm}

\begin{proof}
	Notice that changing coordinate as above locally near each branch point gives
	\begin{align*}
		c_\phi:(\beta,\gamma, \delta)\mapsto
		& \Res\left[ \KS(\beta)\cup\Tr\left( \gamma\cup\delta\right)\right] \\
		= & \Res\left[\frac{\beta_1-\beta_2} {\phi_1-\phi_2} [(\gamma_1-\gamma_2) (\delta_1-\delta_2)+\mbox{ regular terms }]\right],
	\end{align*}
	where the other terms involve $\gamma_i$'s and $\delta_i$'s for $i\geq 3$. The key fact is $W$-equivariance. Choose local coordinate $z$ such that the branch-point is at $z=0$. The reflection in $W$ corresponding to the root producing the branch point acts as multiplication by $-1$ on $\tilde{M}$ near the branch point, so $W$-equivariance implies $\gamma_i(z)=\gamma_i(-z)$ for $i\geq 3$ and similarly for $\delta$. Integrating around the branch-point these terms all cancel themselves out, and what remains is a symmetric formula for the cubic.
\end{proof}


\end{document}